\documentclass{amsart}
\usepackage{booktabs}
\usepackage{makecell}
\usepackage{tabularx}
\usepackage{enumitem,amsmath,amsfonts,amssymb,xcolor}
\usepackage{latexsym,euscript,gensymb}
\newtheorem{lemma}{\bf Lemma}[section]
\newtheorem{nt}[lemma]{\bf Notation}

\newtheorem{thm}[lemma]{\bf Theorem}

\newtheorem{cor}[lemma]{\bf Corollary}

\newtheorem{conj}[lemma]{\bf Conjecture}

\newcommand{\PSL}{{\operatorname{PSL}}}

\newcommand{\PSU}{{\operatorname{PSU}}}

\input xy
\xyoption{all}

\title[]{On Zappa's question in the case of alternating groups}

\author{Ru Zhang}
\address{Ru Zhang. Department of Mathematics, Hubei Minzu University \\ Enshi, Hubei Province,
445000, P. R. China.}
\email{ruzhang2024@hotmail.com}

\author{Rulin Shen}
\address{Rulin Shen. Department of Mathematics, Hubei Minzu University \\ Enshi, Hubei Province,
445000, P. R. China. }
\email{shenrulin@hotmail.com}

\thanks{Project supported by the NSF of China (Grant No. 12161035)}

\date{}



\begin{document}

\setlength{\parskip}{2mm}

\maketitle

\begin{abstract}
In 1962, Guido Zappa asked whether a non-trivial coset of a Sylow $p$-subgroup of a finite group could contain only elements whose orders are powers of $p$.  Marston Conder gives a positive answer to this question in the case of $p=5$. It is known that the smallest group satisfying the conditions of this problem must be a non-abelian simple group.  In this paper, we prove that the smallest group of the Zappa problem could not be an alternating simple group for any prime $p$.
\end{abstract}

\section{Introduction}
Let $G$ be a finite group, and let $P$ be a Sylow $p$-subgroup of $G$ (for some prime $p$
dividing the order of $G$). In a paper in 1962, Guido Zappa \cite{GZ} posed the following two
questions:

\textit{Can a non-trivial coset $Pg$ contain only elements whose orders are powers of $p$? And
in that case, can at least one element of Pg have order $p$}?

A similar question was raised by L.J. Paige:

\textit{Does a non-trivial coset $Pg$ of a Sylow 2-subgroup always contain an element of odd order?}

Paige’s question was answered by John Thompson, who showed in 1967 that for certain primes $q$, such as 53, the group $SL(2,q)$ is a counter-example (see \cite{T}); or in other words, a non-trivial coset of a Sylow 2-subgroup can sometimes contain only elements of even order.
This was taken further in 2014 by Daniel Goldstein and Robert Guralnick, who showed that for every odd prime $p$, there exist infinitely many finite simple groups in which some non-trivial coset $Pg$ of a Sylow $p$-subgroup $P$ contains only elements whose orders are divisible by $p$ (see \cite{GG}). In addition, Daniel Goldstein and Robert Guralnick repeated the first question of Zappa.

In 2017, Marston Conder gives a positive answer to this problem in the case of $p=5$, showing somewhat surprisingly that in a number of finite non-abelian simple groups (including $\PSL(3, 4)$, $\PSU(5, 2)$
and the Janko group $J_3$), some non-trivial coset of a Sylow
5-subgroup (of order 5) contains only elements of order 5 (see \cite{CM}). Moreover, he proved that the smallest finite group with the given property for the prime $p$ must be a non-abelian simple group (see \cite[Lemma 4]{CM}). Kundu and Mishra showed that the smallest finite group is not the alternating group $A_{2p}$ (see \cite{KM}). In this paper, we prove that the smallest group could not be an alternating simple group for all prime $p$. Our main result is as follows.

\begin{thm}
   Let $G$ be an alternating group or a symmetric group, and $P$ be a Sylow $p$-subgroup of $G$. If a coset $P\alpha$  contains only elements whose orders are powers of $p$, then  $\alpha\in P$.
\end{thm}

Simple groups with the given property seem very rare for the prime $p$. Indeed, the forthcoming paper shows that the classical simple group would not occur if $p$ is its characteristic. There exists considerable evidence supporting the validity of the following conjecture posed by Marston Conder, yet it remains open.

\begin{conj}
   Let $S$ be a simple group,  and $P$ be a Sylow $p$-subgroup of $S$. If a nontrivial coset $Pg$  consists entirely of
$p$-elements, then  $|P|=5$.
\end{conj}

\section{Some Notations and Lemmas }
Throughout this paper, the letter $p$ always denotes a prime number.
We denote by $\mathbb{N} = \{0, 1, \dots\}$ the set of natural numbers and by $\mathbb{Z}^+ = \{1, 2, \dots\}$ the set of positive integers. Denoted by $n_p$ the $p$-part of $n$.  Write $S_\Omega$ or $S_n$ ($A_\Omega$ or $A_n$) the symmetric group (alternating group) in the set $\Omega=\{1,2,\cdots, n\}$. Denoted by $id$ the identity mapping, and $\alpha|_{M}$ a permutation restricted in the set $M$. If $\Omega' \subseteq \Omega$, then there exists a natural embedding of $S_{\Omega'}$ in $S_{\Omega}$, so $S_{\Omega'}$ can be regarded as a subgroup of $S_{\Omega}$.  For $k \in \mathbb{N}$ and $0\leq i\leq p-1$, define the set:
    $$\Omega_{k,i} = \{1 + ip^k,\ 2 + ip^k,\ \dots,\ p^k+ip^k\}.$$
    Denoted by $\Omega_k =\Omega_{k,0}$, that is,
    $\Omega_{k}= \{1, 2, \dots, p^k\}.$
We define a mapping $\sigma_k$ in $\Omega_k$ as follows: $\sigma_0=id$ and for $k\geq 1$, $$\sigma_k: x\mapsto x + p^{k-1} \pmod{p^k}.$$

\begin{nt}
 Let  $P_k:=\langle \sigma_0, \sigma_1, \dots, \sigma_k \rangle$ for $k\geq 0$.
 Let $Q_1 = \{ \mathrm{id} \}$, and
    $Q_k := \langle Q_{k-1},\ P_{k-1}^{\sigma_k},\ \dots,\ P_{k-1}^{\sigma_k^{p-1}} \rangle$ for $ k > 1$.
\end{nt}

\begin{lemma}
\label{sigma_k and Omega_k lem}
    For $k \in \mathbb{Z}^+$, the following holds:
    \begin{enumerate}
     \item $\Omega_k = \bigsqcup_{i=0}^{p-1} \Omega_{k-1,i}$.
        \item $\sigma_k \in S_{\Omega_k}$ and $o(\sigma_k) = p$.
        \item $\Omega_{k-1}^{\sigma_k^i} = \Omega_{k-1,i}$ for all $0 \leq i \leq p-1$.
    \end{enumerate}
\end{lemma}

Proof.
Since the elements of $\Omega_k$ are exactly the integers from $1$ to $p^k$, we can partition $\Omega_k$ (increasing order) into $p$ disjoint subsets, each of size $p^{k-1}$. This yields a $p \times p^{k-1}$-matrix by the rule: the first row consists of the integers $1, 2, \dots, p^{k-1}$, and each subsequent row is obtained by adding $p^{k-1}$ to the corresponding elements of the previous row, that is,

\[
    \begin{pmatrix}{}
        1                   & 2 & \dots  & p^{k-1} \\
     1 + p^{k-1}         & 2 + p^{k-1}         & \dots  & 2p^{k-1} \\
      \vdots & \vdots &   &\vdots \\
        1 + (p-1)p^{k-1}  & 2 + (p-1)p^{k-1}  & \dots  & p^k \\
    \end{pmatrix}
\]\\
 where the set of entries of the $i$-th row is just $\Omega_{k-1,i-1}$ for $1 \leq i \leq p$. Therefore, we have
$\Omega_k = \Omega_{k-1} \ \sqcup \ \Omega_{k-1,1} \ \sqcup \ \dots \ \sqcup \ \Omega_{k-1,p-1}.$ This proved (1).

For $1 \leq j \leq p^{k-1}$, by the definition of $\sigma_k$,  we have precisely:
\[
    \sigma_k: j \to j + p^{k-1} \to \dots \to j + (p-1)p^{k-1} \to j \tag{*1}\label{an orb of sigma_k}.
\]
Note that
$\{j, j + p^{k-1}, \dots, j + (p-1)p^{k-1}\}$ is the $j$-th column of the matrix. Thus, the restricted $\sigma_k$ on the $j$-th column is a $p$-cycle. We have $\sigma_k$ is a bijection on $\Omega_k$ and
 $\sigma_k$ is of order $p$, so (2) holds. Moreover, for $1 \leq i \leq p-1 $, \eqref{an orb of sigma_k} implies $\Omega_{k-1,i-1}^{\sigma_k}=\Omega_{k-1,i}$, and so
$\Omega_{k-1}^{\sigma_k^i} = \Omega_{k-1,i}.$ This proved (3). \qed

\begin{lemma}
\label{the p-power of n!}
Let $p$ be a prime and $n \in \mathbb{Z}^+$. Write $n$ to the base $p$: $n = a_0 + a_1^p + \dots + a_mp^m$, where $0 \leq a_i < p$ for each $i$. Let $v_p(n)$ be the exponent of the highest power of $p$ dividing $n!$. Then
$v_p(n) = a_1 + a_2 \frac{p^2-1}{p-1} + \dots + a_m \frac{p^m-1}{p-1}.$
\end{lemma}

Proof. The result follows from \cite[p.49]{D}.\qed

\begin{lemma}
\label{Syl p subg cons}
    Let $k \in \mathbb{N}$ and $P_k = \langle \sigma_0, \sigma_1, \dots, \sigma_k \rangle$ described in Notation 2.1. Then $P_k \in \operatorname{Syl}_p(S_{\Omega_k})$ and 
    $P_k$ $=(P_{k-1}\times P_{k-1}^{\sigma_k}\times P_{k-1}^{\sigma_k^2}\times\cdots\times P_{k-1}^{\sigma_k^{p-1}})\rtimes\langle\sigma_k\rangle$.
\end{lemma}

Proof. By the definition of $P_k$, we have $P_{k+1}=\langle P_k,\sigma_{k+1}\rangle$. The result is clearly true if $k\leq 1$, so proceed by induction on $k$, and assume that $P_k\in Syl_p(S_{\Omega_k})$ is true. We use the following claim to induce the case of $k+1$.

\textbf{Claim:} For $k\geq0$, if $Q\in Syl_p(S_{\Omega_k})$, then $\langle Q, \sigma_{k+1}\rangle\in Syl_p(S_{\Omega_{k+1}})$, and $\langle Q, \sigma_{k+1}\rangle =(Q\times Q^{\sigma_{k+1}}\times Q^{\sigma_{k+1}^2}\times\cdots\times Q^{\sigma_{k+1}^{p-1}})\rtimes\langle\sigma_{k+1}\rangle$.

In fact, by Lemma~\ref{sigma_k and Omega_k lem}, $\Omega_{k}^{\sigma_{k+1}^i} = \Omega_{k,i}$ for $0 \leq i \leq p-1$. Thus, the conjugate subgroup $Q^{\sigma_{k+1}^i}$ is a Sylow $p$-subgroup of $S_{\Omega_{k,i}}$, so is of $S_{\Omega_k}$. Since $\Omega_k = \sqcup_{i=0}^{p-1} \Omega_{k,i}$, we have
$Q Q^{\sigma_{k+1}} \cdots Q^{\sigma_{k+1}^{p-1}} = Q \times Q^{\sigma_{k+1}} \times \cdots \times Q^{\sigma_{k+1}^{p-1}}.$
By Lemma~\ref{sigma_k and Omega_k lem}, $o(\sigma_{k+1}) = p$, so $\sigma_{k+1} \in N_{S_{\Omega_{k+1}}}( \prod_{i=0}^{p-1} Q^{\sigma_{k+1}^i} )$. Furthermore, from $\Omega_k^{\prod_{i=0}^{p-1} Q^{\sigma_{k+1}^i}} = \Omega_k$ and $\Omega_{k}^{\sigma_{k+1}^i} = \Omega_{k,i}$ for $0 \leq i \leq p-1$, it follows that
$\left( \prod_{i=0}^{p-1} Q^{\sigma_{k+1}^i} \right) \cap \langle \sigma_{k+1} \rangle = 1$.
 Hence,
$\langle Q, \sigma_{k+1} \rangle = (Q \times Q^{\sigma_{k+1}} \times \cdots \times Q^{\sigma_{k+1}^{p-1}}) \rtimes \langle \sigma_{k+1} \rangle.$
Thus, $|\langle Q, \sigma_{k+1} \rangle| = |Q|^p \cdot p=(|S_{\Omega_k}|_p)^p\cdot p$. On the other hand,  Lemma~\ref{the p-power of n!} implies
$|\langle Q, \sigma_{k+1} \rangle|=p^{(\frac{p^k-1}{p-1}) \cdot p+1} =|S_{p^{k+1}}|_p
.$
So $\langle Q, \sigma_{k+1} \rangle \in \operatorname{Syl}_p(S_{\Omega_{k+1}})$.
\qed\medskip

Using Lemma~\ref{sigma_k and Omega_k lem} and Lemma~\ref{Syl p subg cons}, the following corollary can be obtained directly.

\begin{cor}
\label{syl-p subg of S_Omeg(k,i)}
Let $k \in \mathbb{Z}^+$ and $0 \leq i<p$. Then $P_{k-1}^{\sigma_k^i} \in Syl_p(S_{\Omega_{k-1,i}})$.
\end{cor}

\begin{lemma}
\label{Q_k}
   Let $k \in \mathbb{Z}^+, k > 1$.
    Then $Q_k \leq P_k$ and $Q_k \in \operatorname{Syl}_p(S_{\Omega_k \setminus \{1\}})$.
\end{lemma}

Proof. According to the definition of $P_k$,  $P_{k-1}^{\sigma_k^i} \leq P_k$ for $k \in \mathbb{Z}^+$. It follows that $Q_k \leq P_k$. Thus, by Lemma~\ref{Syl p subg cons}, we have
\[
Q_k = Q_{k-1} \times P_{k-1}^{\sigma_k} \times \dots \times
P_{k-1}^{\sigma_k^{p-1}}. \tag{*2}\label{Q_k cons}
\]
Since $|\Omega_k| = p^k$ and $|\Omega_k \setminus \{1\}| = p^k - 1 = (p-1)(p^{k-1} + \dots + p + 1)$, Lemma~\ref{the p-power of n!} implies that
\[
|S_{\Omega_k}|_p = p^{\frac{p^k - 1}{p - 1}}
\tag{*3}\label{|S_p^k|_p}
\]
and
\[
|S_{\Omega_k \setminus \{1\}}|_p = p^{(p-1)\left(1 +
\frac{p^2-1}{p-1} + \dots + \frac{p^{k-1}-1}{p-1}\right)}.
\tag{*4}\label{|no 1|_p}
\]
We use induction to prove $Q_k \in \operatorname{Syl}_p(S_{\Omega_k \setminus \{1\}})$. When $k = 1$, we have $Q_1 = \{\mathrm{id}\}$ and $|S_{\Omega_1 \setminus \{1\}}|_p = 1$, so $Q_1 \in \operatorname{Syl}_p(S_{\Omega_1 \setminus \{1\}})$ is true.  Since $P_{k-1}^{\sigma_k^i} \leq S_{\Omega_{k-1,i}}$ and the elements of $S_{\Omega_{k-1,i}}$ fix $1$ for $0 < i < p$, by the equality \eqref{Q_k cons}, we have
$Q_k \leq S_{\Omega_k \setminus \{1\}}$
and
$|Q_k| = |Q_{k-1}| \cdot |P_{k-1}|^{p-1}$.  By the hypothesis that $Q_{k-1} \in \operatorname{Syl}_p(S_{\Omega_{k-1} \setminus \{1\}})$, we have
$|Q_k| = |S_{\Omega_{k-1} \setminus \{1\}}|_p \cdot (|S_{\Omega_{k-1}}|_p)^{p-1}.$
The equalities \eqref{|S_p^k|_p} and \eqref{|no 1|_p} give
$|Q_k| = p^{(p-1)\left(1 + \frac{p^2-1}{p-1} + \dots + \frac{p^{k-2}-1}{p-1}\right)} \cdot p^{p^{k-1} - 1} = p^{(p-1)\left(1 + \frac{p^2-1}{p-1} + \dots + \frac{p^{k-1}-1}{p-1}\right)} = |S_{\Omega_k \setminus \{1\}}|_p.$
Thus, $Q_k \in \operatorname{Syl}_p(S_{\Omega_k \setminus \{1\}})$. \qed

\begin{lemma}
\label{syl-p subg of S_n}
Let $\Omega' \subseteq \Omega$ with $|\Omega'|=p^k$. Write $n$ in base $p$:
$n = a_0 + a_1p +\dots + a_kp^k +\dots + a_mp^m,$
where $a_i\in \mathbb{N}$, $0 \leq a_i <p$,  and $a_k \neq 0$. If $P \in \operatorname{Syl_p(S_{\Omega'})}$ and $Q \in \operatorname{Syl_p(S_{\Omega \setminus \Omega'})}$, then $P\times Q \in \operatorname{Syl_p(S_\Omega)}$.
\end{lemma}
Proof. It follows from Lemma~\ref{the p-power of n!} that the Sylow $p$-subgroups of $S_n$ have order $p^{v_p(n)}$, where
$v_p(n) = a_1 + a_2 \frac{p^2-1}{p-1} + \dots + a_k \frac{p^k-1}{p-1} \dots + a_m \frac{p^m-1}{p-1}.
$
Let $n' = n - p^k$. Then $|\Omega \setminus \Omega'| = n'$ and $n' = a_0 + a_1p +\dots + (a_k-1)p^k +\dots + a_mp^m,$
thus,
$v_p(n') = a_1 + a_2 \frac{p^2-1}{p-1} + \dots + (a_k-1) \frac{p^k-1}{p-1} \dots + a_m \frac{p^m-1}{p-1}.$
Therefore, $v_p(n) = v_p(n')+v_p(p^k)$, and hence $|P\times Q|=|S_n|_p$.\qed

\begin{lemma}
\label{a prop of sly subg}
Let $G$ be a group and $P \in \operatorname{Syl}_p(G)$. Then the following statements are equivalent:

\begin{enumerate}
    \item For all $\alpha \in G$, if every element of the coset $P\alpha$ is a $p$-element, then $\alpha \in P$.
    \item For all $Q \in \operatorname{Syl}_p(G)$ and all $\alpha \in G$, if every element of the coset $Q\alpha$ is a $p$-element, then $\alpha \in Q$.
\end{enumerate}
\end{lemma}

Proof. (2) $\Rightarrow$ (1) is obvious. It suffices to prove (1) $\Rightarrow$ (2).
Take any $\alpha \in G$ and $Q \in \operatorname{Syl}_p(G)$. By Sylow's Theorem, there exists $g \in G$ such that $P = Q^g$. Let $\beta = \alpha^g$. Then every element in $P\beta$ is a $p$-element, so by (1) we have $\beta \in P$, and hence $\alpha \in Q$.\qed

\begin{lemma}
\label{number}
    Let $s>1$ and $a>1$ be integers, and let $t$ and $t_i$ be natural numbers for $1\leq i\leq s$. If $a^{t_1}+a^{t_1}+\cdots +a^{t_s}=a^t$, then $s\geq a$, and the equality $s=a$ holds if and only if $t_1=t_2=\cdots =t_s=t-1$.
\end{lemma}

Proof. Set $t_M:=\max\{t_1,t_2,\cdots, t_s\}$. Since $s>1$, $a^{t_i}<a^t$ for all $1\leq i\leq s$, and so $t_M<t$. It follows that  $a^t=a^{t_1}+a^{t_2}+\cdots+a^{t_s}\leq s\cdot a^{t_M}\leq s\cdot a^{t-1}$, so we have $s\geq a$. Moreover, if $s=a$, then $a^t=a^{t_1}+a^{t_2}+\cdots+a^{t_s}\leq a\cdot a^{t-1}=a^t$, and so $a^{t_i}=a^{t_M}=a^{t-1}$ for all $i$.\qed

\section{Symmetric Groups}

In this section, we discuss the symmetric groups.

\begin{lemma}
\label{act trans}
    Let $P \in \operatorname{Syl}_p(S_{p^k})$. Then $P$ acts transitively on $\Omega_k$.
\end{lemma}
Proof. Let $\alpha = (1, 2, \dots, p^k)$ be a $p^k$-cycle. Then $\alpha$ is a $p$-element. By Sylow's Theorem, there exists $x \in P$ such that $\alpha^x \in P$. Since $\langle \alpha \rangle$ acts transitively on $\Omega_k$, so does $\langle \alpha^x \rangle$ on $\Omega_k$. Therefore, $P$ acts transitively on $\Omega_k$.\qed

\begin{lemma}
\label{orb lem}
Let $\alpha \in S_n$ and $T$ be an orbit of $\langle \alpha \rangle$. Let $a \in T$ and $l \in \mathbb{Z}^+$ with $l \mid |T|$. Then:
\begin{enumerate}
    \item Let $T_i = \{a^{\alpha^{i+kl}} \colon 0 \leq k < \frac{|T|}{l} \}$, $0 \leq i < l$. Then $T_i $ is an orbit of $\langle\alpha^l|_T\rangle$ for $0\leq i<l$ and $ T = \sqcup_{i=0}^{l-1} T_i$.
    \item Let $T' = \{ a, a^{\alpha^{t_1}}, \dots, a^{\alpha^{t_s}} \}$, where $0 < t_1 < t_2 < \dots < t_s < |T|$. If $T'$ is an orbit of $\langle \alpha^l \rangle$, then $s = \frac{|T|}{l} - 1$, and $t_i = i l$ for $1 \leq i \leq s$.
\end{enumerate}
\end{lemma}
Proof. (1).  Let $a^{\alpha^i \langle \alpha^l \rangle}=\{a^{\alpha^{i}\alpha^j}|\alpha^j\in\langle \alpha^l\rangle\}.$ We claim that $T_i=a^{\alpha^i \langle \alpha^l \rangle}$. In fact,  $T_i \subseteq a^{\alpha^i \langle \alpha^l \rangle} $ is obvious. Conversely, taking any $ b\in a^{\alpha^i \langle \alpha^l \rangle}$, there exists $h\in \mathbb{N}$, such that $b=a^{\alpha^{i+hl}}$. By the division algorithm, there exists $h_0 ,r\in\mathbb{N}$ with $\ 0\leq r<\frac{|T|}{l} $ such that $h=h_0\cdot\frac{|T|}{l}+r$. Thus,
$b= a^{\alpha^{i+(h_0\cdot\frac{|T|}{l}+r)l}} =a^{\alpha^{i+rl}} \in T_i,$ and so $a^{\alpha^i \langle \alpha^l \rangle} \subseteq T_i$. So $T_i$ is just an orbit of $\langle\alpha^l|_T\rangle$, and then $ T = \sqcup_{i=0}^{l-1} T_i$.

(2). According to (1), $T_i$ is an orbit of $\langle \alpha^l\rangle$, for $0\leq i <l $. Since $a\in T_0 \cap T' \neq \emptyset$, we have  $T'=T_0$. Hence, $s=\frac{|T'|}{l}-1 = \frac{|T|}{l}-1$. Moreover, since $0 < t_1 < t_2 < \dots < t_s < |T|$, we have $t_i= il$, for $ 1\leq i \leq s$.\qed

\begin{lemma}
\label{mult-orbt lem}
Let $H \subseteq \Omega$, $\alpha, \beta \in S_n$, with $\beta|_{\Omega\setminus H}=\mathrm{id}$. Let $T_\alpha$ be an orbit of $\langle \alpha \rangle$ and $T_{\beta\alpha}$ be an orbit of $\langle \beta\alpha \rangle$. If $H \subseteq T_{\beta\alpha}$ and $H \cap T_\alpha \neq \emptyset$, then $T_\alpha \subseteq T_{\beta\alpha}$.
\end{lemma}

Proof. Take any $a \in T_\alpha$. Since $H \cap T_\alpha \neq \emptyset$ and the action of $\langle \alpha \rangle$ on $T_\alpha$ is transitive, for all $b\in H \cap T_\alpha$, there exist $k_0 \in \mathbb{Z}^+$ such that
$b^{\alpha^{k_0}} = a.$  We claim that  there exist $b\in H\cap T_\alpha$ and $k\in Z^+$ such that  $b^{\alpha^{k}}=a$ and  $b^{\alpha^j}\notin H$  for all $0<j<k$. In fact,  fix $b\in H \cap T_\alpha$. We may assume that there exists $i \,(0<i<k_0)$ such that  $b^{\alpha^i} \in H$, without loss of generality, the $i$ is the largest satisfied $b^{\alpha^i} \in H$, then $b^{\alpha^j}\notin H$ for $i<j<k_0$. Let $b_1=b^{\alpha^i}$ and $k=k_0-i$.  Then $b_1^{\alpha^{k}}=a$  and $b_1^{\alpha^j}\notin H$  for all $0<j<k$.
Since $\beta|_{\Omega \setminus H} = \mathrm{id}$, it follows that, for all $0 < j< k$, $\beta$ fixes $b^{\alpha^j}$, that is,  $$(b^{\alpha^j})^\beta=b^{\alpha^j}.\eqno (*5)$$
Also, because $\beta|_H\in S_H$, there exists $b_0 \in H$ such that $b_0^\beta = b$. Thus (*5) implies $b_0^{(\beta\alpha)^k}=b_0^{\beta\alpha(\beta\alpha)^{k-1}}$ $=b^{\alpha(\beta\alpha)^{k-1}}$
$=((b^{\alpha})^\beta)^{\alpha(\beta\alpha)^{k-2}}$$=(b^{\alpha})^{\alpha(\beta\alpha)^{k-2}}$$=(b^{\alpha^2})^{(\beta\alpha)^{k-2}}$ $=\cdots=b^{\alpha^k}$, and so $b_0^{(\beta\alpha)^k}=a$.
 Therefore, $b_0$ and $a$ belong to the same orbit of $\langle \beta\alpha \rangle$. Since $b_0 \in H\subseteq T_{\beta\alpha}$, it follows that $a \in T_{\beta\alpha}$. By our hypothesis $a\in T_\alpha$, we have $T_\alpha \subseteq T_{\beta\alpha}$.\qed

\begin{lemma}
\label{first thm}
    Let $k  \in \mathbb{Z^+}$, $\Omega'\subseteq \Omega$ with $|\Omega'|=p^k$, and $P \in  Syl_p(S_{\Omega'})$. Let $\alpha \in S_{\Omega}$ and $T$ be an $\langle \alpha\rangle$-orbit.
    If $T\cap \Omega'\neq \emptyset$ and all elements in the coset $P\alpha$ are $p$-elements, then $|T\cap \Omega'|=1$ or $|T\cap \Omega'|\geq p$ and the equality holds if and only if $T\cap \Omega'$ is a
     $\langle \alpha^{\frac{|T|}{p}}\rangle$-orbit.
\end{lemma}

Proof. Since $T$ is an orbit of $\langle \alpha \rangle$, taking $a \in T$, we have
$\alpha|_T = (a, a^\alpha, \dots, a^{\alpha^{|T|-1}}).$
Let $s = |T \cap \Omega'| > 0$. Since $T \cap \Omega' \subseteq T$, we can assume that $T \cap \Omega' = \{a^{\alpha^{u_1}}, \dots, a^{\alpha^{u_s}}\}$ and $0\leq u_1 < u_2 < \dots < u_s<|T|$.
 Denote by $a_i = a^{\alpha^{u_i}}$ for $1 \leq i \leq s$, and $l_i = u_{i+1} - u_i$ for $1 \leq i \leq s-1$, and $l_s = u_1 - u_s + |T|$. It follows that $T \cap \Omega' = \{a_1, \dots, a_s\}$ and we may rewrite
$$\alpha|_T = (a_1, a_1^\alpha, \dots, a_1^{\alpha^{l_1-1}},  \dots, a_s, a_s^\alpha, \dots, a_s^{\alpha^{l_s-1}})
, \eqno (*6)$$
where $ a_i^{\alpha^k} \notin T \cap \Omega'$  for $ 0 < k < l_i$ and $ 1 \leq i \leq s$, and $a_i^{\alpha^{l_i}} = a_{i+1} $ for $1\leq i<s$ and $ a_s^{\alpha^{l_s}} = a_{1}$. For convenience,  let $a_{s+1}=a_1$ in the following.
Let $K_i = \{a_i^\alpha, \dots, a_i^{\alpha^{l_i-1}}, a_{i+1}\}$ for $1 \leq i \leq s$, then $|K_i| = l_i$, and $T = \sqcup_{i=1}^s K_i$. Since  $|\Omega'|=|\Omega_k|=p^k$, $S_{\Omega_k}\cong S_{\Omega'}$.  By Lemma \ref{act trans}, there exists $\beta_i\in P$ such that $a_{i+1}^{\beta_i} = a_i$ for $1 \leq i \leq s$. Also, $\beta_i|_{\Omega \backslash \Omega'} = \mathrm{id}$ implies that $\beta_i$ fixes $a_i^{\alpha^k}$ for all $0 < k < l_i $. Therefore,
$a_{i+1}^{(\beta_i \alpha)^k} = a_i^{\alpha^k} $ for $0 < k \leq l_i.$ In particular, (*6) implies $\alpha_i^{\alpha^{l_i}}=\alpha_{i+1}$.
So, $K_i$ is an orbit of $\langle \beta_i \alpha \rangle$, and $(\beta_i \alpha)|_{K_i} = (a_{i+1}, a_i^\alpha, a_i^{\alpha^2}, \dots, a_i^{\alpha^{l_i-1}}).$
Since all elements of $P\alpha
$ are $p$-elements, $\beta_i \alpha$ is a $p$-element. Thus, the size $|K_i|=l_i=p^{t_i}$ for some $t_i \geq 0$. Similarly, since $\alpha$ is a $p$-element, there exists $t \geq 0$ such that $|T| = p^t$. Moreover, $T = \sqcup_{i=1}^s K_i$, so we have $|T| = l_1 + \dots + l_s$, i.e.,
$p^t = p^{t_1} + \dots + p^{t_s}.$
If $s > 1$, Lemma \ref{number} implies $s \geq p$,  and if $s = p$,  then $t_i = t - 1$ for all $1 \leq i \leq s$. Thus
$a_{i+1} = a_i^{\alpha^{l_i}} = a_i^{\alpha^{p^{t-1}}},$
and hence
$T \cap \Omega' = \{a_1^{\alpha^{i p^{t-1}}} | 0 \leq i < p\}.$
By Lemma \ref{orb lem}, $T \cap \Omega'$ is an orbit of $\langle \alpha^{p^{t-1}} \rangle$, so is $\langle\alpha^\frac{|T|}{p} \rangle$.
Conversely, if $T \cap \Omega'$ is an orbit of $\langle\alpha^\frac{|T|}{p} \rangle$,  then item (2) of Lemma \ref{orb lem} leads to $|k\cap\Omega'|=p$.
\qed

\begin{cor}
\label{first cor}
Let $\alpha \in S_n$, $\beta = (b_1, \dots, b_p) \in S_n$. If all elements in $\langle \beta \rangle \alpha$ are  $p$-elements, then either
\begin{enumerate}
    \item $b_1, \dots, b_p$ belong to $p$ distinct orbits of $\langle \alpha \rangle$, or
    \item there exists an orbit $T$ of $\langle \alpha \rangle$ such that $\{b_1, \dots, b_p\} \subseteq T$, and $\{b_1, \dots, b_p\}$ is an orbit of $\langle \alpha^{|T|/p} \rangle$.
\end{enumerate}
\end{cor}

Proof. Let $\Omega'=\{b_1,\dots,b_p\}$. The results follow from Lemma \ref{first thm}.\qed

\begin{thm}
\label{orbt cons}
 Let $k\in \mathbb{N}$, $\alpha \in S_n$ and $T$ be an orbit of $\langle \alpha \rangle$. Let  $\Omega' \subseteq T$ and $|\Omega'| = p^k$, and let $P \in \text{Syl}_p(S_{\Omega'})$. If all elements of the coset $P\alpha$ are  $p$-elements, then $\Omega'$ is an orbit of $\langle \alpha^{{|T|}/{p^k}} \rangle$.
\end{thm}
Proof. We proceed by induction on $k$. If $k = 0$, the conclusion holds obviously, so we assume that the conclusion holds for case $k - 1$. Since $|\Omega'| = |\Omega_k|=p^k$, there exists a bijection $\tau \colon \Omega' \to \Omega_k$. Since $P \in \text{Syl}_p(S_{\Omega'})$, we have $P^{\tau} \in \text{Syl}_p(S_{p^k})$. Also, since $P_k \in \text{Syl}_p(S_{p^k})$, there exists $h \in S_{p^k}$ such that $P_k = P^{\tau h}$. Let $\alpha_0 = \alpha^{\tau h}$ and $H = T^{\tau h}$. Since $\Omega' \subseteq T$ and $T$ is an orbit of $\langle\alpha\rangle$, it follows that $\Omega_k \subseteq H$ and $H$ is an orbit of $\langle\alpha_0\rangle$. Furthermore,  all elements in the coset $P\alpha$ are  $p$-elements, so is in the coset $P_k \alpha_0$.
Let $|H| = p^t$ and  $l = p^{t - k + 1}$. Take $a \in H$, then $H = \{ a^{\alpha_0^i} | 0 \leq i < |H| \}$. Denote by $H_i = \{ a^{\alpha_0^{i + jl}} | 0 \leq j < {|H|}/{l} \}$ for all $0 \leq i < l$. By Lemma \ref{orb lem}, $H_0, $ $\dots,$ $ H_{l-1}$ are all orbits of $\langle \alpha_0^l|_H \rangle$.
Also, let $c_{ij} = a^{\alpha_0^{i + jl}}$ for $0 \leq i < l$, $0 \leq j < {|H|}/{l}$. Then we have $H_i = \{ c_{i,j} | 0 \leq j < {|H|}/{l} \}$, and
$$\alpha_0|_H = (c_{0,0}\ c_{1,0}\dots c_{l-1,0}\ c_{0,1}\ c_{1,1} \dots c_{l-1,1}\dots c_{0,{|H|}/{l} - 1}\ c_{1,{|H|}/{l} - 1} \dots c_{l-1,{|H|}/{l} - 1}). \eqno (*7)$$
 By Lemma~\ref{Syl p subg cons},
$P_{k-1} \leq P_{k}$ and $\sigma_k\in P_k$, so $P_{k-1}^{\sigma_k^s}\leq P_k$  for all {{$0 \leq s < p$. }}
Since all elements in $P_k \alpha_0$ are $p$-elements, it follows that all elements in $P_{k-1}^{\sigma_k^s} \alpha_0$ are $p$-elements. Furthermore, Corollary \ref{syl-p subg of S_Omeg(k,i)} implies $P_{k-1}^{\sigma_k^s} \in \operatorname{Syl}_p(S_{\Omega_{k-1,s}})$.  The induction hypothesis gives that $\Omega_{k-1,s}$ is an orbit of $\langle \alpha_0^l|_H \rangle$. Without loss of generality,  we assume $\Omega_{k-1,s} = H_{i_s}$ and $0 = i_0 < i_1 < \dots < i_{p-1} < l$. Hence, Lemma~\ref{sigma_k and Omega_k lem} yields
$\Omega_k = \sqcup_{s=0}^{p-1} H_{i_s}.$
For convenience, let $c_{0,|H|/l}=c_{0,0}$. Now,  for $0 \leq j\leq |H|/l-1$, we denote $K_j = \{ c_{1,j}, \dots, c_{l-1,j}, c_{0,j+1} \}$. Then $K_j \cap H_0 = \{ c_{0,j+1} \}$ and $K_j \cap H_i = \{ c_{i,j} \}$ for $0 < i < l$.
Thus, $$K_j \cap \Omega_k = \sqcup_{s=0}^{p-1} (K_j \cap H_{i_s}) = \{c_{0,j+1}, c_{i_1,j}, \dots, c_{i_{p-1},j}\}.\eqno (*8)$$
 Since $c_{0,0}, c_{0,1} \in H_0 = \Omega_{k-1}$, by Lemma \ref{act trans}, there exists $g \in P_{k-1}$ such that $c_{0,1}^{g} = c_{0,0}$. Since $g|_{\Omega \setminus \Omega_{k-1}} = \mathrm{id}$ and $\{ c_{i,0} | 0 \leq i < l \} \cap \Omega_{k-1} = \{ c_{0,0} \}$, it follows that $g$ fixes $c_{i,0}$ for $0 < i < l$. Hence, (*7) implies that
 $c_{0,1}^{(g \alpha_0)^i} = c_{i,0} $ for $1 \leq i < l$
and $c_{0,1}^{(g \alpha_0)^l} = c_{0,1}.$
Thus, $K_0$ is an orbit of $\langle g \alpha_0 \rangle$, and $$(g \alpha_0)|_{K_0} = (c_{0,1}, c_{1,0}, \dots, c_{l-1,0}).\eqno (*9)$$
Since every element in $P_k (g \alpha_0)$ is a $p$-element and $|K_0 \cap \Omega_k| = p$, Lemma~\ref{first thm} implies that $K_0 \cap \Omega_k$ is an orbit of $\langle (g \alpha_0)^{{|K_0|}/{p}} \rangle$. Note that (*8) and (*9) give $K_0 \cap \Omega_k = \{ c_{0,1}, c_{i_1,0}, \dots, c_{i_{p-1}, 0} \} = \{ c_{0,1}^{ (g \alpha_0)^s } | s = 0, i_1, \dots, i_{p-1}\}.
$
By Lemma~\ref{orb lem} we have $i_s =s{|K_0|}/{p}= s{l}/{p}$  for $0 \leq s< p$. Therefore,
$\Omega_k = \sqcup_{s=0}^{p-1} \{ a^{\alpha_0^{i_s + j l}} | 0 \leq j < {|H|}/{l} \} = \sqcup_{s=0}^{p-1} \{ a^{\alpha_0^{(s + j p) l / p}} | 0 \leq j < {|H|}/{l} \}$ $= \{ a^{\alpha_0^{i p^{t-k}}} | 0 \leq i \leq p^k - 1 \}.$
As Lemma~\ref{orb lem}, $\Omega_k$ is an orbit of $\langle \alpha_0^{{|H|}/{p^k}} \rangle$, and hence $\Omega'$ is an orbit of $\langle \alpha^{{|T|}/{p^k} }\rangle$.\qed

\begin{lemma}
\label{non empt lem}
Let $p > 2$, $1 < r < p$ and $k \in \mathbb{Z}^+$. Let $\alpha \in S_\Omega$,  $T$ be an orbit of $\langle \alpha \rangle$, and $\Omega_k \subseteq\Omega $. Suppose
$\sqcup_{i=1}^{r} \Omega_{k-1,t_i} \subseteq T$,
where $t_1,t_2,\dots,t_r$ are distinct elements of $\{1,\dots,p\}$. If every element in $P_k \alpha$ is a $p$-element, then
$
T \cap \left( \Omega_k \setminus \sqcup_{i=1}^{r} \Omega_{k-1,t_i} \right) \neq \emptyset.
$
\end{lemma}

Proof. By Lemma~\ref{Syl p subg cons}, $P_k = (P_{k-1} \times P_{k-1}^{\sigma_k} \times \dots \times P_{k-1}^{\sigma_k^{p-1}}) \rtimes \langle \sigma_k \rangle$.
Denoted by  $P_{k-1,t} = P_{k-1}^{\sigma_k^t}$ for $0 \leq t < p$. Then
$P_k = (P_{k-1,0} \times \dots \times P_{k-1,p-1}) \rtimes \langle \sigma_k \rangle,$
and
$P_{k-1,t} \in \operatorname{Syl}_p(S_{\Omega_{k-1,t}}).$
Let $l = {|T|}/{p^{k-1}}$. Take $a \in T$ and define
$T_s = \{ a^{\alpha^{s + jl}} | 0 \leq j < {|T|}/{l}\}$ for $0 \leq s < l.$
By Lemma~\ref{orb lem}, $T_0, \dots, T_{l-1}$ are all orbits of $\langle \alpha^l|_T \rangle$.
Denote by $c_{s,j} = a^{\alpha^{s + jl}}$ for $0 \leq s < l$ and $0 \leq j < {|T|}/{l}$. Then $T_s = \{ c_{s,j} | 0 \leq j < {|T|}/{l}\}$ for $0 \leq s < l$, and
\[
\alpha|_T = (c_{0,0}, \dots, c_{l-1,0},\ c_{0,1}, \dots, c_{l-1,1},\ \dots,\ c_{0,{|T|}/{l}-1}, \dots, c_{l-1,{|T|}/{l}-1}).\eqno (*10)
\]
Since $P_{k-1,t_i} \leq P_k$ for $1\leq i\leq r $, every element in $P_{k-1,t_i} \alpha$ is a $p$-element. Also, since $\Omega_{k-1,t_i} \subseteq T$, Theorem~\ref{orbt cons} implies that $\Omega_{k-1,t_i}$ is an orbit of $\langle \alpha^l|_T \rangle$. Therefore, there exists $0 \leq u_i < l$ such that $\Omega_{k-1,t_i} = T_{u_i}$. Without loss of generality, assume $0 = u_1 < u_2 < \dots < u_r < l$. Denote by $c_{0,|T|/l}:=c_{0,0}$.  Define
$K_j = \{c_{1,j}, \dots, c_{l-1,j}, c_{0,j+1}\}$ for $ 0 \leq j < {|T|}/{l}$.
Then
$T = \sqcup_{j=0}^{{|T|}/{l}-1} K_j,$
and
\[
K_j \cap \left( \sqcup_{i=1}^r \Omega_{k-1,t_i} \right) = K_j \cap \left( \sqcup_{i=1}^r T_{u_i} \right) = \{c_{u_2,j}, \dots, c_{u_r,j}, c_{0,j+1}\}.\eqno (*11)
\]
Now we choose $j=0$, the equality (*11) implies \[
|K_0 \cap \Omega_k| \geq \left| K_0 \cap \left( \sqcup_{i=1}^r \Omega_{k-1,t_i} \right) \right| = r > 1. \eqno (*12)
\]
Since $c_{0,0}, c_{0,1} \in T_0 = \Omega_{k-1,t_1}$, by Lemma~\ref{act trans} there exists $g \in P_{k-1,t_1}$ such that $c_{0,1}^g = c_{0,0}$. Since $c_{s,0} \in T \setminus T_0$ for $0 < s < l$, the element $g$ fixes $c_{s,0}$. Therefore,
(*10) follows $c_{0,1}^{(g\alpha)^s} = c_{s,0} $ for $ 1 \leq s < l,$
and
$c_{0,1}^{(g\alpha)^l} = c_{0,1}.$
Thus,
$(g\alpha)|_{K_0}$ $= (c_{0,1}, c_{1,0}, \dots, c_{l-1,0}).$ In particular, $K_0$ is an orbit of $\langle g\alpha \rangle$.
Since all elements in $P_k(g\alpha)$ are $p$-elements and
$|K_0 \cap \Omega_k| > 1,
$
Lemma~\ref{first thm} gives $|K_0 \cap \Omega_k| \geq p$. By (*12)
$|K_0 \cap \Omega_k| = r + \left| K_0 \cap \left( \Omega_k \setminus \sqcup_{i=1}^r \Omega_{k-1,t_i} \right) \right|$
and $r < p$, so
$\left| K_0 \cap \left( \Omega_k \setminus \sqcup_{i=1}^r \Omega_{k-1,t_i} \right) \right| > 0.$
Moreover, since $K_0 \subseteq T$, we conclude that
$\left| T \cap \left( \Omega_k \setminus \sqcup_{i=1}^r \Omega_{k-1,t_i} \right) \right| > 0.$ \qed

\begin{lemma}
\label{capacity of orbit}
    Let $\Omega'\subseteq \Omega$ and $|\Omega'|=p^k$, let $P\in Syl_p(S_{\Omega'})$ and $\alpha\in S_{\Omega}$. If every element of the coset $P\alpha$ is a $p$-element, then there exists $\gamma\in P$ such that $\Omega'$ is a subset of some $\langle \gamma\alpha\rangle$-orbit.
\end{lemma}
Proof.  If $k = 0$, then $|\Omega'|=1$. Taking $\gamma = \mathrm{id}$, the result is obviously true, so we assume $k>0$ and work by induction on $k$.

Since $|\Omega'| = |\Omega_k|$, there exists a bijection $\tau \colon \Omega' \to \Omega_k$, and so $P^\tau \in \operatorname{Syl}_p(S_{\Omega_k})$. Note that $P_k \in \operatorname{Syl}_p(S_{\Omega_k})$. Then there exists $h \in S_{\Omega_k}$ such that $P_k = P^{\tau h}$. Let $\alpha_0 := \alpha^{\tau h}$.  It follows that  every element in $P_k \alpha_0$ is a $p$-element.
 For $0 \leq i < p$, we set $P_{k-1,i} := P_{k-1}^{\sigma_k^i}$. By Lemma~\ref{sigma_k and Omega_k lem} and Lemma~\ref{syl-p subg of S_Omeg(k,i)},  $\Omega_k = \sqcup_{i=0}^{p-1} \Omega_{k-1,i}$ and $P_{k-1,i} \in \operatorname{Syl}_p(S_{\Omega_{k-1,i}})$. Moreover, Lemma~\ref{Syl p subg cons} gives $P_{k-1,i}\leq P_k$. Let $\delta \in P_k$, we have $P_{k-1,i}g\alpha_0 \subseteq P_k\alpha_0$, so every element in $P_{k-1,i}(\delta \alpha_0)$ is a $p$-element.
The inductive hypothesis for the case of $k-1$ guarantees that

(*H): \textit{For $g \in P_k$, there exists $\gamma \in P_{k-1,i}$ such that $\Omega_{k-1,i}$ is contained in an orbit of $\langle \gamma g\alpha_0 \rangle$.} \medskip
\\
Next, we use the $(k-1)$-hypothesis (*H) to prove the following two claims.

\textbf{Claim 1:} \textit{Let $t$ be a positive number less than $p$ and $\gamma \in P_k$. Let $T$ be an orbit of $\langle \gamma \alpha_0 \rangle$ and $i_1,\cdots, i_t$ are different natural numbers less than $p$. If $\sqcup_{j=1}^{t} \Omega_{k-1,i_j} \subseteq T $, then there exist $\gamma'\in P_k$ and an orbit $T'$ of $\langle\gamma'\alpha_0\rangle$ such that $\sqcup_{j=1}^{t} \Omega_{k-1,i_j} \subseteq T' $ and $T' \cap \left( \Omega_k \setminus \sqcup_{j=1}^{t} \Omega_{k-1,i_j} \right) \neq \emptyset$. }

In fact, if $t >1$, let $\gamma'= \gamma $ and $T'=T$, then $\sqcup_{j=1}^{t}\Omega_{k-1,i_j}\subseteq T'$. Furthermore, since $t>1$ and all elements in $P_k(\gamma'\alpha_0)$ are $p$-elements, Lemma~\ref{non empt lem} yields $T' \cap \left( \Omega_k \setminus \sqcup_{j=1}^{t} \Omega_{k-1,i_j} \right) \neq \emptyset$. {If $t=1$, we may assume }that $T_1 \cap \left( \Omega_k\setminus \Omega_{k-1,i_1}\right) = \emptyset$. Let $j \neq i_1$ and $0 \leq j < p$, and choose  points $a \in \Omega_{k-1,i_1}$, $b \in \Omega_{k-1,j}$. By Lemma~\ref{act trans}, there exists $g \in P_k$ such that $b^g = a$. Since $T$ is an orbit of $\langle \gamma \alpha_0 \rangle$ and $\Omega_{k-1,i_1} \subseteq T$, we have $a^{\gamma \alpha_0} \in T$. Also, since $\Omega_{k-1,i_1} \subseteq T$, there exists $s \in \mathbb{Z}^+$ such that $a^{(\gamma \alpha_0)^s} \in \Omega_{k-1,i_1}$, and $a^{(\gamma \alpha_0)^i} \notin \Omega_{k-1,i_1}$ for $0 < i < s$. Hence,  $T_1 \cap $ $\left( \Omega_k\setminus \Omega_{k-1,i_1}\right) = \emptyset$ implies that $a^{(\gamma \alpha_0)^i} \notin \Omega_k$ for $0 < i < s$. Moreover, since $g|_{\Omega \setminus \Omega_k} = \mathrm{id}$,
$a^{(\gamma \alpha_0)^i g} = a^{(\gamma \alpha_0)^i}, \quad 0 < i < s.$
Thus, $b^{(g \gamma \alpha_0)^i} = a^{(\gamma \alpha_0)^i}, \quad 0 < i \leq s.$
In particular, $b^{(g \gamma \alpha_0)^s} = a^{(\gamma \alpha_0)^s} \in \Omega_{k-1,i_1}$. Therefore, $b$ and $a^{(\gamma \alpha_0)^s}$ belong to a same orbit of $\langle g \gamma \alpha_0 \rangle$.
Since  $g \gamma \in P_k$ , by the hypothesis (*H) there exist $\delta\in P_{k-1,i_1}$ and an orbit $T'$ of $\langle \delta (g \gamma \alpha_0) \rangle$ such that $\Omega_{k-1,i_1} \subseteq T'$. Since $b$ and $a^{(\gamma\alpha_0)^s}$ belong to the same orbit of $\langle g \gamma \alpha_0 \rangle$, and $a^{(\gamma \alpha_0)^s} \in \Omega_{k-1,i_1}$, Lemma~\ref{mult-orbt lem} yields  $b \in T'$. Thus, $T'\cap(\Omega_k\setminus\Omega_{k-1,i_1})\neq\emptyset$. Let $\gamma'=\delta g\gamma_1$, then $\gamma'\in P_k$ and $T'$ is an orbit of $\langle \gamma'\alpha_0 \rangle$.

 \textbf{Claim 2:} \textit{For $ 1\leq t \leq p $, there exist distinct numbers $ 0\leq i_1, \dots, i_t <p$ and $\gamma_t \in P_k $
and an orbit $T_t$ of $\langle\gamma_t\alpha_0\rangle$ such that
$\sqcup_{j=1}^{t}\Omega_{k-1,i_j}\subseteq T_t$. }

In fact, if $t=1$, it is just hypothesis (*H), then the result is already true.  We assume that the result holds for the case of $t-1$, that is, there exist $t-1$ different integers $i_1,\cdots, i_{t-1}$,  $\gamma_{t-1} \in P_k$  and an orbit $T_{t-1}$ of $\langle \gamma_{t-1}\alpha_0 \rangle$ such that $\sqcup_{j=1}^{t-1}\Omega_{k-1,i_j}\subseteq T_{t-1}$, and we will prove it holds for the case of $t$. By Claim 1,
there exist $\gamma' \in P_k$  and an orbit $T'$ of $\langle \gamma'\alpha_0 \rangle$ such that $\sqcup_{j=1}^{t-1}\Omega_{k-1,i_j}\subseteq T'$ and $T' \cap \left( \Omega_k \setminus \sqcup_{j=1}^{t-1} \Omega_{k-1,i_j} \right) \neq \emptyset$.
So, there exists $i_t \in \{1, \dots, p\} \setminus \{i_j | 1 \leq j \leq t-1\}$ such that $T' \cap \Omega_{k-1,i_{t}} \neq \emptyset$. Since $\gamma' \in P_k $, by the induction hypothesis (*13) there exist $\gamma_t' \in P_{k-1,i_t}$ and an orbit $T_t$ of $\langle \gamma_t' (\gamma' \alpha_0) \rangle$ such that $\Omega_{k-1,i_t} \subseteq T_t$. Furthermore,  since $\Omega_{k-1,i_t} \cap T' \neq \emptyset$ and $\gamma_t'|_{\Omega \setminus \Omega_{k-1,i_t}} = \mathrm{id}$, Lemma~\ref{mult-orbt lem} implies that $T' \subseteq T_t$, hence $\sqcup_{j=1}^{t} \Omega_{k-1,i_j} \subseteq T_t$. Let $\gamma_t = \gamma_t' \gamma'$, then $\gamma_t \in P_k$ and $T_t$ be an orbit of $\langle \gamma_t \alpha_0 \rangle$. Therefore, we have proved Claim 2.

In particular, when $t=p$, there exist $\gamma_p$ and an orbit $T_p$ of $\langle \gamma_t\alpha_0\rangle$,  such that $\Omega_k =\sqcup_{j=1}^{p}\Omega_{k-1,i_j} \subseteq T_p$, this holds for the case of $k$. As required.
 \qed

\begin{thm}
\label{p-elements thm}
    Let $\alpha \in S_n$ and $P \in \operatorname{Syl}_p(S_n)$. If all elements in the coset $P\alpha$ are  $p$-elements, then $\alpha \in P$.
\end{thm}

Proof. We note that if $p\nmid |S_n|$, then the trivial subgroup is a Sylow $p$-subgroup.   In the following, we fix the prime $p$ and $\Omega:=\{1,2,\cdots, n\}$, $\Omega_m:=\{1,2,\cdots, p^m\}$. The $P_m$ and $Q_m$ are the same as in Notation 2.1, and the $S_\Omega$ means $S_n$.

If $n = 1$ and 2, the results hold, obviously. We assume that $n>2$ and proceed by induction on $n$.
Write $n$ in base $p$:
$n = a_0 + a_1p + \dots + a_mp^m,$
where $a_m \neq 0$ and $0 \leq a_i < p$ for all $i$.
Let $n' = n - p^m$. Then $|\Omega \setminus \Omega_m| = n'$. Choose $Q \in \operatorname{Syl}_p(S_{\Omega \setminus \Omega_m})$. Lemma~\ref{syl-p subg of S_n} implies $P_m \times Q \in \operatorname{Syl}_p(S_n)$. By Lemma~\ref{a prop of sly subg}, we may assume that $P = P_m \times Q$.

Case 1. $n' = 0$, i.e, $n=p^m$.  Then $\Omega = \Omega_m$ and $P = P_m$. Since $1,1^{\alpha^{-1}} \in \Omega_m $, by Lemma \ref{act trans} there exists $g \in P$ such that $1^g = 1^{\alpha^{-1}}$, so $g\alpha \in S_{\Omega_m \setminus \{1\}}$.   Lemma~\ref{Q_k} leads to $Q_m \leq P$,  so all elements in $Q_m(g\alpha)$ are $p$-elements. Also, by Lemma~\ref{Q_k}  $Q_m \in \operatorname{Syl}_p(S_{\Omega_m \setminus \{1\}})$, and then the induction hypothesis gives $g\alpha \in Q_m \leq P$. Thus $\alpha \in P$.

Case 2. $n'\neq 0 $. Since $P=P_m \times Q$, it is easy to see that $\exp(P)=\exp(P_m)=p^m $. As every element in the coset $P\alpha$ is a p-element, every element in $P_m\alpha$ is also a $p$-element. By  Lemma~\ref{capacity of orbit}, there exist $\beta \in P_m$ and an orbit $T$ of $\langle \beta\alpha \rangle $, such that $\Omega_m \subseteq T$. Moreover, since $\beta\alpha$ is a $p$-element, $o(\beta\alpha)\leq p^m$. Thus,
$p^m=|\Omega_m| \leq |T| \leq o(\beta\alpha) \leq p^m,$
which implies that $|T|=
|\Omega_m| = p^m.$
Hence $T= \Omega_m$. Set $\alpha_1=(\beta\alpha)|_{\Omega_m}$ and $\alpha_2=(\beta\alpha)|_{\Omega \setminus \Omega_m}$. Then $\alpha_1 \in S_{\Omega_m}, \alpha_2 \in S_{\Omega \setminus \Omega_m} $ and $\beta\alpha=\alpha_1\alpha_2$. Since every element in $P(\beta\alpha)$ is a $p$-element, every element in $P_m\alpha_1$ and $Q\alpha_2$ is also a $p$-element. As $p^m<n$ and $ n'< n$, the induction hypothesis yields $\alpha_1\in P_m, \alpha_2\in Q$. Therefore, $\beta\alpha \in P$, and then $\alpha\in P$.\qed

\section{Alternating Groups}

In this section, we discuss the case of alternating groups. The notation is the same as in Section 3.

\begin{lemma}
\label{'p-cyc' mov 'p pt'}
Let $c_1,\dots,c_p$ be distinct elements in $\Omega$ and $P \in \operatorname{Syl}_p(S_n)$. Then there exists a $p$-cycle $\beta\in S_n$ such that at least one point in $c_1\dots,c_p$ is moved by $\beta$.
\end{lemma}
Proof. We first prove the result

(*R):
\textit{if $k > 0$, $c \in \Omega_k$, and $P \in \operatorname{Syl}_p(S_{p^k})$, then there exists a $p$-cycle $\beta \in P$ such that $c$ is moved by $\beta$.}

In fact, if $k = 1$, $P$ is generated by a $p$-cycle. Lemma~\ref{act trans} implies that $P$ acts transitively on $\Omega_1$. Thus, there exists $\beta \in P$ such that $c$ is moved by $\beta$. We assume that $k>1$ and proceed by induction on $k$.
Now, suppose that the conclusion holds for $k-1$. As Sylow subgroups are conjugate, we may assume $P = P_k$. By Lemmas~\ref{sigma_k and Omega_k lem} and~\ref{Syl p subg cons}, we have $\Omega_k = \sqcup_{i=0}^{p-1} \Omega_{k-1,i}$ and $P_{k-1}^{\sigma_k^i}\leq P_k$ for all $0 \leq i < p$.  Since $c\in \Omega_k$,
there exists some $i$ such that $c \in \Omega_{k-1,i}$. Moreover, Lemma~\ref{syl-p subg of S_Omeg(k,i)} implies that $P_{k-1}^{\sigma_k^i} \in \operatorname{Syl}_p(S_{\Omega_{k-1,i}})$. By the induction hypothesis, there exists a $p$-cycle $\beta \in P_{k-1}^{\sigma_k^i}$ such that $c$ is moved by $\beta$.

Now, we prove the lemma by induction. Write $n$ in base $p$:
$n = a_0 + a_1p + \dots + a_mp^m,$
where $a_m \neq 0$ and $0 \leq a_i < p$ for all $i$.

If $|\Omega| = n = p$, the conclusion holds by the result (*R) proved above. Now suppose that $n>p$ and work by induction on $n$.   Let $n' = n - p^m$ and $Q \in \operatorname{Syl}_p(S_{\Omega \setminus \Omega_m})$. Certainly, $|\Omega\setminus \Omega_m|=n'.$
By Lemma~\ref{syl-p subg of S_n}, we have $P_m \times Q \in \operatorname{Syl}_p(S_n)$. Again, by conjugacy, we may assume $P = P_m \times Q$.
If some $c_i \in \Omega_m$, then the result follows from (*R). If $c_i \in \Omega \setminus \Omega_m$ for all $1 \leq i \leq p$, then the induction hypothesis on $n'$ implies that there exists a $p$-cycle $\beta \in Q$ such that $\beta$ moves some point of the set $\{c_1,\cdots, c_p\}$.
\qed

Next, we consider the Sylow 2-subgroups of alternating groups by $S_4$. Recall $\Omega_1=\{1,2\}$, $\Omega_{1,1}=\{3,4\}$, $\Omega_2=\{1,2,3,4\}$ and $P_2=\langle\sigma_1, \sigma_2 \rangle$ as in Notation 2.1.
Easily compute $P_2=\{\mathrm{id},(12),(34),(12)(34),(13)(24),(14)(23),(1324),(1423)\}$, that is just a Sylow 2-subgroup of $S_4$.

\begin{lemma}
\label{p_2 prop}
Let $p=2$, $\alpha \in A_n$, and $P\in\operatorname{Syl}_2(S_n)$ with $P_2\leq P$. If every element in $(P \cap A_n)\alpha$ is a $2$-element, then every element in $P_2\alpha$ is a $2$-element.
\end{lemma}

Proof. Denote by $i^{\langle\alpha\rangle}$ an $\langle\alpha\rangle$-orbit containing $i$.  Let $T = \cup_{i=1}^{4} i^{\langle\alpha\rangle}$. Then $\alpha|_T \in S_T$, $\alpha|_{\Omega \setminus T} \in S_{\Omega \setminus T}$, and $\alpha = \alpha|_T \cdot \alpha|_{\Omega \setminus T}$. Taking any $g \in P_2$, we have $g\alpha = (g\alpha|_T) \cdot (\alpha|_{\Omega \setminus T})$. Since $\alpha|_{\Omega \setminus T}$ is a $2$-element, $g\alpha$ is a $2$-element if and only if $g\alpha|_T$ is a $2$-element. Therefore, it suffices to prove that every element of $P_2(\alpha|_T)$ is a $2$-element. It is obvious that $P_2 \cap A_n=\{ \mathrm{id},(12)(34),(13)(24),(14)(23) \}$ acts transitively on $\Omega_2$. Note that $T$ is a union of at most 4 distinct  $\langle\alpha\rangle$-orbits.

\noindent \textbf{Case 1.} Each  element of $\Omega_2$ lies in distinct $\langle \alpha \rangle$-orbits, i.e, $T$ has 4 orbits.

Let $i\in \Omega_2$ and $|i^{\langle \alpha \rangle}| = l_i$.  Since $\alpha$ is a $2$-element,  we may assume that $l_i = 2^{t_i}$ for some $t_i \in \mathbb{N}$.
Note that $i^{\langle \alpha \rangle} = \{i, i^\alpha, \dots, i^{\alpha^{l_i-1}}\}$.  Since $(12)(34)$ fixes any point in the set $i^{\langle\alpha\rangle}\setminus\{i\}$ for all $i\in\Omega_2$. We have
$(12)(34)\alpha|_T=
(1^\alpha,$ $\dots,$ $1^{\alpha^{l_1-1}},1,$ $2^\alpha,\dots,2^{\alpha^{l_2-1}},2)
(3^\alpha,\dots,3^{\alpha^{l_3-1}},3,4^\alpha,\dots,4^{\alpha^{l_4-1}},4)$.
Hence, $1^{\langle \alpha \rangle} \sqcup 2^{\langle \alpha \rangle}$ and $3^{\langle \alpha \rangle} \sqcup 4^{\langle \alpha \rangle}$ are orbits of $\langle(12)(34)\alpha\rangle$. Since $(12)(34)\alpha\in (P_2\cap A_n)\alpha$, $(12)(34)\alpha$ is a $2$-element, and so  $l_1 + l_2 = 2^t$ and $l_3 + l_4 = 2^{t'}$ for some  $t, t' \in \mathbb{Z^+}$. It follows that $2^{t_1} + 2^{t_2} = 2^t$ and $2^{t_3} + 2^{t_4} = 2^{t'}$. By  Lemma~\ref{number},  $l_1 = l_2$ and $l_3 = l_4$.
By similar discussion of $(13)(24)\alpha$, we obtain $l_1 = l_3$ and $l_2 = l_4$. Hence, $l_1 = l_2 = l_3 = l_4$, which is a power of $2$, say $2^t$. Note that $P_2=\{\mathrm{id},(12),$ $(34),(12)(34),(13)(24),(14)(23),(1324),(1423)\}$. It is straightforward to verify that for every $g \in P_2$, the orbit size of $\langle g\alpha|_T \rangle$ is $2^t$, $2^{t+1}$, or $2^{t+2}$, and hence all elements of $P_2(\alpha|_T)$ are $2$-elements.

\noindent \textbf{Case 2.} Exactly three points of $\Omega_2$ belong to a single $\langle\alpha\rangle$-orbit.

Let  $T_1$ and $T_2$ be two  $\langle \alpha \rangle$-orbits with $|T_1 \cap \Omega_2| = 3$ and $|T_2 \cap \Omega_2| = 1$. Now we set $T_1 \cap \Omega_2 = \{i_1, i_2, i_3\}$ and $T_2 \cap \Omega_2 = \{i_4\}$.
Then $\Omega_2=\{i_1,i_2,i_3,i_4\}$. There exist positive integers $k_1, k_2$ with  $0 < k_1, k_2 < |T_1|$ such that
$i_1^{\alpha^{k_1}} = i_2 \quad \text{and} \quad i_1^{\alpha^{k_2}} = i_3.$
Without loss of generality, assume $k_2 > k_1$. Let $l_1 = k_1$, $l_2 = k_2 - k_1$, $l_3 = |T_1| - k_2$, and $l_4 = |T_2|$. Then we have
\[
\alpha|_T = (i_1\ i_1^\alpha, \dots, i_1^{\alpha^{l_1-1}}\ i_2\ i_2^\alpha\ \dots\ i_2^{\alpha^{l_2-1}}\ i_3\ i_3^\alpha\ \dots\ i_3^{\alpha^{l_3-1}})(i_4\ i_4^\alpha\ \dots\ i_4^{\alpha^{l_4-1}}).
\]
It follows that
\[
(i_1 i_2)(i_3 i_4)\alpha|_T = (i_4\ i_3^\alpha \dots i_3^{\alpha^{l_3-1}}\ i_1\ i_2^\alpha \dots i_2^{\alpha^{l_2-1}}\ i_3\ i_4^\alpha \dots i_4^{\alpha^{l_4-1}})(i_2\ i_1^\alpha \dots i_1^{\alpha^{l_1-1}}),
\]
\[
(i_1i_3)(i_2i_4)\alpha|_T = (i_3\ i_1^\alpha \dots i_1^{\alpha^{l_1-1}}\ i_2\ i_4^\alpha\dots i_4^{\alpha^{l_4-1}}\ i_4\ i_2^\alpha \dots i_2^{\alpha^{l_2-1}})(i_1\ i_3^\alpha \dots i_3^{\alpha^{l_3-1}}),
\]
\[
(i_1i_4)(i_2i_3)\alpha|_T = (i_1\ i_4^\alpha \dots i_4^{\alpha^{l_4-1}}\ i_4\ i_1^\alpha \dots i_1^{\alpha^{l_1-1}}\ i_2\ i_3^\alpha \dots i_3^{\alpha^{l_3-1}})(i_3\ i_2^\alpha \dots i_2^{\alpha^{l_2-1}}).
\]
Since every element in $(P_2 \cap A_n)\alpha $ is a $2$-element, the orbit length of each element in $(P_2 \cap A_n)\alpha $ is a power of $2$. From the cycle structure of above four permutations, it follows the equations
\begin{align}
l_1 + l_2 + l_3 = 2^{s_4}, ~~~~~~&~~~~~l_4 = 2^{t_4}, \tag{E1}\label{equ-7} \\
l_2 + l_3 + l_4 = 2^{s_1}, ~~~~~~&~~~~~ l_1 = 2^{t_1}, \tag{E2}\label{equ-8} \\
l_1 + l_2 + l_4 = 2^{s_3}, ~~~~~~&~~~~~ l_3 = 2^{t_3}, \tag{E3}\label{equ-9} \\
l_1 + l_3 + l_4 = 2^{s_2}, ~~~~~~&~~~~~l_2 = 2^{t_2}, \tag{E4}\label{equ-10}
\end{align}
for some $s_i \in \mathbb{Z}^+$, $t_i \in \mathbb{N}$ with $ 1 \leq i \leq 4$. From the equations \eqref{equ-7} \eqref{equ-8}, we have
$l_1-l_4=2^{s_4}-2^{s_1}=2^{t_1}-2^{t_4},$ which implies that $2^{s_1}(2^{s_4-s_1}-1)=2^{t_4}(2^{t_1-t_4}-1).$
If $s_4 \neq s_1$, then $s_1 = t_4$. So the equation \eqref{equ-8} becomes $l_2 + l_3 + l_4 = l_4$, hence $l_2 = l_3 = 0$, a contradiction. If $s_4 = s_1$, similarly, from the equations \eqref{equ-7} \eqref{equ-9}, and \eqref{equ-7} \eqref{equ-10}, we obtain $s_4 = s_3$ and $s_4 = s_2$, respectively. Thus $s_1=s_2=s_3=s_4$. Computing the sums of the left and the right from (E1) to (E4) respectively, we get
$3(l_1+l_2+l_3+l_4)=4\cdot2^{s_4} ,$
which is a contradiction. Therefore, Case 2 is impossible.

\noindent \textbf{Case 3.} The four elements of $\Omega_2$ lie in three distinct $\langle\alpha\rangle$-orbits, i.e, $T$ has 3 orbits.

Let $T_1$, $T_2$, and $T_3$ be  three distinct $\langle \alpha \rangle$-orbits, and let $T_1 \cap \Omega_2 = \{i_1, i_2\}$, $T_2 \cap \Omega_2 = \{i_3\}$, and $T_3 \cap \Omega_2 = \{i_4\}$.  As in Case 2, we may get that
\[
\alpha|_T = (i_1\ i_1^\alpha \dots i_1^{\alpha^{l_1-1}}\ i_2\
i_2^\alpha \dots i_2^{\alpha^{l_2-1}})(\ i_3\ i_3^\alpha \dots
i_3^{\alpha^{l_3-1}})(i_4\ i_4^\alpha \dots
i_4^{\alpha^{l_4-1}}).\tag{*13}\label{equ-11}
\]
where $l_1+l_2=|T_1|$, $l_3=|T_2|$, and $l_4=|T_3|$. Then
\[(i_1i_2)(i_3i_4)\alpha|_T = (i_1 i_2^\alpha \dots i_2^{\alpha^{l_2-1}})(i_2 i_1^\alpha \dots i_1^{\alpha^{l_1-1}})(i_3 i_4^\alpha \dots i_4^{\alpha^{l_4-1}} i_4\ i_3^\alpha \dots i_3^{\alpha^{l_3-1}}),\tag{*14}\label{equ-12} \]
$$(i_1i_3)(i_2i_4)\alpha|_T = (i_3\ i_1^\alpha \dots i_1^{\alpha^{l_1-1}}\ i_2\ i_4^\alpha \dots i_4^{\alpha^{l_4-1}}\ i_4\ i_2^\alpha \dots i_2^{\alpha^{l_2-1}}\ i_1\ i_3^\alpha \dots i_3^{\alpha^{l_3-1}}). ~~~$$
Since every element in $(P_2\cap A_n)\alpha$ is a $2$-element, we obtain the equations $l_1 + l_2 = 2^{s_1}, $ $ l_3 + l_4 = 2^{s_2},$ $
l_i = 2^{t_i}  (1 \leq i \leq 4)$ and $
l_1 + l_2 + l_3 + l_4 = 2^s$
for some $s_1,s_2,s\in \mathbb{Z}^+$ and $t_i\in \mathbb{N}$ with $ 1\leq i\leq4$.
By Lemma~\ref{number}, we obtain $l_1=l_2=l_3=l_4.$ Denote by $l:=l_1$, and denote by $K_j := \{i_j^{\alpha^k} | 1\leq k\leq l\}$ for $j=1,2$ and   $K_3:=T_2$,  $K_4:=T_3.$
Note that $(T_1\cap\Omega_2)\sqcup(T_2\cap\Omega_2)\sqcup(T_3\cap\Omega_2)=\Omega_2=\Omega_1\sqcup \Omega_{1,1}$.  We claim that

(*C)$:\ \ \  T_1\cap\Omega_2=\Omega_1$ or $T_1\cap\Omega_2=\Omega_{1,1}$.

Otherwise, $|T_1\cap \Omega_1|=1$ and $|T_1\cap\Omega_{1,1}|=1.$ Without loss of generality, we may assume that $T_1\cap \Omega_1=\{i_1\}$, $T_1\cap \Omega_{1,1}=\{i_2\}$, $T_2\cap\Omega_1=\{i_3\}$, and $T_3\cap\Omega_{1,1}=\{i_4\}$. Then $(i_1i_3)=(1\ 2)\in P_2$ and $(i_2i_4)=(3\ 4)\in P_2.$ Since $\alpha$ is an  even permutation and $\alpha|_T$ is an odd permutation, it follows that $\alpha|_{\Omega\setminus T} \neq \mathrm{id}$, and then $|\Omega\setminus T|>1$. By Lemma~\ref{'p-cyc' mov 'p pt'}, there exists a $2$-cycle $\beta \in P$ such that $\beta$ moves some point in $\Omega\setminus T$. Let $\beta=(j_1j_2)$, where $j_1 \in \Omega\setminus T$. Thus,
$$(i_1i_3)(j_1j_2),\ (i_2i_4)(j_1j_2) \in P\cap A_n.\eqno (*15)$$
Since $T = \sqcup_{i=1}^{4}K_j$, we have $j_2\notin T$ or $j_2$ belongs to some $K_j$, and so $j_2\notin K_1\sqcup K_2$ or $j_2\notin K_3\sqcup K_4$. We now consider two cases accordingly.

 \textbf{(I).}  $j_2\notin K_1\sqcup K_2$.
Then there exists $u\in\{3, 4\}$ such that $j_2\notin K_u$.
If $u=3$, we have
$(i_1 i_3)(j_1j_2)\alpha|_{\bigsqcup_{i=1}^{3}K_i}=(i_3 i_1^\alpha\dots i_1^{\alpha^{l-1}} i_2 i_2^\alpha\dots i_2^{\alpha^{l-1}} i_1 i_3^\alpha\dots i_3^{\alpha^{l-1}}).$
Thus, $\sqcup_{i=1}^{3}K_i$ is an orbit of $\langle(i_1i_3)(j_1j_2)\alpha\rangle$ and $|(\sqcup_{i=1}^{3}K_i) \cap\Omega_2|=3$. Note that all elements in $(P\cap A_n)(i_1i_3)(j_1j_2)\alpha$ are $2$-elements, as $(i_1i_3)(j_1j_2)\in P\cap A_n$ from (*15). It is impossible by Case 2. Similarly, if $u=4$, then
$(i_2 i_4)(j_1j_2)\alpha|_{K_1\sqcup K_2\sqcup K_4}=(i_1 i_1^\alpha\dots i_1^{\alpha^{l-1}}$ $ i_2 i_4^\alpha\dots i_4^{\alpha^{l-1}} i_4 i_2^\alpha\dots i_2^{\alpha^{l-1}}),$
and so $K_1\sqcup K_2\sqcup K_4$ is an orbit of $\langle(i_2i_4)(j_1j_2)\alpha\rangle$, which is also impossible by Case 2.

\textbf{(II).}  $j_2\notin K_3\sqcup K_4.$
 Let $\alpha_0=(i_1,i_2)(i_3,i_4)\alpha$. Since $(i_1i_2)(i_3i_4)\in P_2\cap A_n\subseteq P\cap A_n$, all elements in $(P\cap A_n)\alpha_0$ are $2$-elements.  By the cycle structure of $\alpha_0$ in  \eqref{equ-12}, we have  $K_1$,$K_2$, and $K_3\sqcup K_4$ are orbits of $\langle\alpha_0\rangle$. Note that $|(K_3\sqcup K_4)\cap\Omega_1|=1$, $|(K_3\sqcup K_4)\cap\Omega_{1,1}|=1$, and $j_2\notin K_3\sqcup K_4$. Let $T_1':=K_3\sqcup K_4$, $K_1':=K_3$ and $K_2':=K_4$. The conditions that $|T_1'\cap \Omega_1|=|T_1'\cap\Omega_{1,1}|=1$ and $j_2 \notin K_1'\sqcup K_2'$ force us back to case (I), which is impossible. Thus, (*C) is proved.

Therefore, (*C) implies that $T_1\cap \Omega_2=\{i_1,i_2\}=\Omega_1$ or $\Omega_{1,1}$, and a direct computation by (*13) shows that every element in the coset $P_2(\alpha|_T)$ is a $2$-element.

 \noindent \textbf{Case 4.} The four elements of $\Omega_2$ form two pairs, each pair lying in one $\langle\alpha\rangle$-orbit.

 Let $T_1$ and $T_2$ be two distinct orbits of $\langle \alpha \rangle$, and  let $T_1 \cap \Omega_2=\{i_1,i_2\}$ and $T_2\cap \Omega_2=\{i_3,i_4\}$. As previous cases,  assume that
$\alpha|_T = (i_1 i_1^\alpha \dots i_1^{\alpha^{l_1-1}} i_2 i_2^\alpha \dots i_2^{\alpha^{l_2-1}})(\ i_3 i_3^\alpha \dots$ $i_3^{\alpha^{l_3-1}} i_4 i_4^\alpha \dots i_4^{\alpha^{l_4-1}}),$
where $l_1+l_2=|T_1|$ and $l_3+l_4=|T_2|$. Then
$(i_1i_2)(i_3 i_4)\alpha|_T = (i_2 i_1^\alpha\dots i_1^{\alpha^{l_1-1}})(i_1 i_2^\alpha\dots i_2^{\alpha^{l_2-1}})(i_4\ i_3^\alpha\dots i_3^{\alpha^{l_3-1}})(i_3 i_4^\alpha \dots i_4^{\alpha^{l_4-1}}).$ Note that all elements in $(P \cap A_n)(i_1i_2)(i_3 i_4)\alpha$ are $2$-elements and  each  element of $\Omega_2$ lies in different orbits of $\langle(i_1i_2)(i_3 i_4)\alpha \rangle$. So, this case is, in fact, Case 1, and thus the result holds.

\noindent \textbf{Case 5.} All elements of $\Omega_2$ lie in a single  $\langle\alpha\rangle$-orbit.

Let $T$ be an orbit of $\langle\alpha\rangle$ and $\Omega_2 \subseteq T$, and let
$\alpha|_T = (i_1 i_1^\alpha \dots i_1^{\alpha^{l_1-1}} i_2 i_2^\alpha \dots i_2^{\alpha^{l_2-1}}$ $ i_3 i_3^\alpha \dots i_3^{\alpha^{l_3-1}} i_4 i_4^\alpha \dots i_4^{\alpha^{l_4-1}}),$
where $\{i_1,i_2,i_3,i_4\}=\Omega_2$. Then,
$
(i_1i_4)(i_2i_3)\alpha|_T = (i_4 i_1^\alpha \dots i_1^{\alpha^{l_1-1}} i_2 i_3^\alpha \dots i_3^{\alpha^{l_3-1}})(i_3 i_2^\alpha \dots i_2^{\alpha^{l_2-1}})(i_1 i_4^\alpha \dots i_4^{\alpha^{l_4-1}}).
$ Thus, this case is, in fact, Case 3, and so the result holds.
\qed

\begin{lemma}
\label{odd-pro}
Let $\alpha \in A_n$ and $P \in \operatorname{Syl}_p(S_n)$. If every element of $(P \cap A_n)\alpha$ is a $p$-element, then every element of $P\alpha$ is a $p$-element.
\end{lemma}
Proof. If $p>2$, then $P \cap A_n = P$. Hence, the result follows from Theorem \ref{p-elements thm}, immediately.
Next, we consider the case $p=2$. It suffices to prove that every element in $(P \setminus (P\cap A_n))\alpha$ is a $2$-element. By Lemma \ref{a prop of sly subg}, we may assume $P_2\leq P$. Take any odd permutation $g \in P$. Since $(1\ 2) \in P_2$, we have $(1\ 2)g \in P \cap A_n$. Thus, every element in $(P \cap A_n)(1\ 2)g\alpha$ is a $2$-element. Lemme~\ref{p_2 prop} implies that all elements in $P_2(1,2)g\alpha $ are $2$-elements, and so $g\alpha$ is a $2$-element. Since $g$ is an arbitrary odd permutation in $P$,  every element in $(P \setminus (P\cap A_n))\alpha$ is a $2$-element.
\qed

Finally, we apply Theorem~\ref{p-elements thm} and Lemma~\ref{odd-pro}, the following result can be obtained at once.

\begin{thm}
Let $\alpha\in A_n$ and $P \in \operatorname{Syl}_p(A_n)$. If all elements in $P\alpha$ are  $p$-elements, then $\alpha \in P$.
\end{thm}

\end{document}